\documentclass[preprint,12pt]{elsarticle}

\usepackage{amssymb}
\usepackage{amsmath}
\usepackage{mathtools}
\usepackage{amsthm}
\usepackage{color}
\usepackage{multirow}
\usepackage{lscape}

\usepackage{subfigure}
\usepackage{float}

\usepackage[svgnames]{xcolor}
\usepackage{soul}

\usepackage{algorithm2e}  

\usepackage{tikz}
\tikzset{every label/.style={font=\footnotesize,inner sep=1pt}}
\newcommand{\stencilpt}[4][]{\node[circle,fill,draw,inner sep=1.5pt,label={above left:#4},#1] at (#2) (#3) {}}

\newtheorem{theorem}{Theorem}

\newcommand{\kmin}{\kappa_{\rm min}}
\newcommand{\kmax}{\kappa_{\rm max}}
\journal{Journal of Computational Physics}

\begin{document}

\begin{frontmatter}



\title{On the equivalence between the Scheduled Relaxation Jacobi method
  and Richardson's non-stationary method}


\author{J.E.~Adsuara\corref{cor1}\fnref{label1}} 
\author{I.~Cordero-Carri\'on\corref{cor1}\fnref{label2}}
\author{P.~Cerd\'a-Dur\'an\corref{cor1}\fnref{label1}}
\author{V.~Mewes\corref{cor1}\fnref{label1}} 
\author{M.A.~Aloy\corref{cor1}\fnref{label1}}
\cortext[cor1]{jose.adsuara@uv.es, isabel.cordero@uv.es, pablo.cerda@uv.es, vassilios.mewes@uv.es, miguel.a.aloy@uv.es}
\address[label1]{Departamento de Astronom\'{\i}a y Astrof\'{\i}sica, Universidad de Valencia, E-46100, Burjassot, Spain.}
\address[label2]{Departamento de Matem\'atica Aplicada, Universidad de Valencia, E-46100, Burjassot, Spain.}

 \begin{abstract}
   The Scheduled Relaxation Jacobi (SRJ) method is an extension of the
   classical Jacobi iterative method to solve linear systems of
   equations ($Au=b$) associated with elliptic problems. It inherits
   its robustness and accelerates its convergence rate computing a set
   of $P$ relaxation factors that result from a minimization
   problem. In a typical SRJ scheme, the former set of factors is
   employed in cycles of $M$ consecutive iterations until a prescribed
   tolerance is reached. We present the analytic form for the optimal
   set of relaxation factors for the case in which all of them are
   strictly different, and find that the resulting algorithm is
   equivalent to a non-stationary generalized Richardson's method
   where the matrix of the system of equations is preconditioned
   multiplying it by $D=\rm{diag}(A)$. Our method to estimate the
   weights has the advantage that the explicit computation of the
   maximum and minimum eigenvalues of the matrix $A$ (or the
   corresponding iteration matrix of the underlying weighted Jacobi
   scheme) is replaced by the (much easier) calculation of the maximum
   and minimum frequencies derived from a von Neumann analysis of the
   continuous elliptic operator. This set of weights is also the
   optimal one for the general problem, resulting in the fastest
   convergence of all possible SRJ schemes for a given grid
   structure. The amplification factor of the method can be found
   analytically and allows for the exact estimation of the number of
   iterations needed to achieve a desired tolerance. We also show that
   with the set of weights computed for the optimal SRJ scheme for a
   fixed cycle size it is possible to estimate numerically the optimal
   value of the parameter $\omega$ in the Successive Overrelaxation
   (SOR) method in some cases. Finally, we demonstrate with 
   practical examples that our method also works very well for
   Poisson-like problems in which a high-order discretization of the Laplacian
   operator is employed (e.g., a $9-$ or $17-$points discretization). This
   is of interest since the former discretizations do not yield
   consistently ordered $A$ matrices and, hence, the theory of Young
   cannot be used to predict the optimal value of the SOR
   parameter. Furthermore, the optimal SRJ schemes
   deduced here are advantageous over existing SOR implementations for
   high-order discretizations of the Laplacian operator in as much as
   they do not need to resort to multi-coloring schemes for their
   parallel implementation. 
\end{abstract}


\end{frontmatter}


\section{Introduction}
\label{sec1}

The Jacobi method \cite{Jacobi1845} is an iterative method to solve
systems of linear equations. Due to its simplicity and its
convergence properties it is a popular choice as preconditioner, in
particular when solving elliptic partial differential equations.
However, its slow rate of convergence, compared to other iterative
methods (e.g. Gauss-Seidel, SOR, Conjugate gradient, GMRES), makes it
a poor choice to solve linear systems. The scheduled relaxation Jacobi
method \cite{Yang2014}, SRJ hereafter, is an extension of the
classical Jacobi method, which increases the rate of convergence in
the case of linear problems that arise in the finite difference
  discretization of elliptic equations. It consists of executing a
series of weighted Jacobi steps with carefully chosen values for the
weights in the sequence. Indeed, the SRJ method can be expressed
  for a linear system, $Au=b$, as
\begin{equation}
u^{n+1} = u^n + \omega_{n} D^{-1}( b - A u^n),
\label{eq:SRJ}
\end{equation}
where $D$ is the diagonal of the matrix $A$.  If we consider a set of $P$
different relaxation factors, $\omega _{n}, \,\, n=1,\dots,P$, such that
$\omega _{n} > \omega _{n+1}$ and we apply each relaxation factor $q_n$
times, the {\it total
  amplification factor} after $M:=\sum_{n=1}^P q_n$ iterations is
\begin{equation}
G_M(\kappa) = \prod_{n=1}^{P} (1 - \omega _{n} \kappa)^{q_n}, \label{eq:gm}
\end{equation}
which is an estimation of the reduction of the residual during one
cycle (M iterations). In the former expression $\kappa$ is a function
of the wave-numbers obtained from a von Neumann analysis of the system
of linear equations resulting from the discretization of the original
elliptical problem by finite differences (for more details see
\cite{Yang2014,Adsuetal15}). Yang \& Mittal \cite{Yang2014} argued
that, for a fixed number $P$ of different weights, there is an optimal
choice of the weights $\omega _{n}$ and repetition numbers $q_n$ that
minimizes the maximum {\it per-iteration amplification factor},
$\Gamma(\kappa) = |G(\kappa)|^{1/M}$, in the interval $\kappa \in
[\kmin,\kmax]$ and therefore also the number of iterations needed for
convergence. The boundaries of the interval in $\kappa$ correspond to
the minimum and the maximum weight numbers allowed by the
discretization mesh and boundary conditions used to solve the elliptic
problem under consideration.

In the aforementioned paper, \cite{Yang2014} computed numerically
the optimal weights for $P\le 5$ and Adsuara et al.\,\cite{Adsuetal15}
extended the calculations up to $P=15$. 
The main properties of the SRJ, obtained by  \cite{Yang2014} and confirmed by \cite{Adsuetal15}, are the following:

\begin{enumerate}
\item Within the range of $P$ studied, increasing the number of weights $P$ improves the rate of convergence.

\item The resulting SRJ schemes converge significantly faster than the classical Jacobi method by factors exceeding $100$ in the 
methods presented by  \cite{Yang2014} and  $\sim 1000$ in those
presented by \cite{Adsuetal15}. Increasing grid sizes, i.e.  decreasing $\kmin$, results in larger acceleration factors.

\item The optimal schemes found use each of the weights multiple times, resulting in a total number of iterations $M$ per
cycle significantly larger than $P$, e.g. for $P=2$, \cite{Yang2014} found an optimal scheme with $M=16$ for the
smallest grid size they considered ($N=16$), while for larger grids $M$ notably increases (e.g., $M=1173$ for $N=1024$).
\end{enumerate}

The optimization procedure outlined by \cite{Yang2014} has a caveat
though.  Even if the amplification factor were to reduce monotonically
by increasing $P$, for sufficiently high values of $P$, the number of
iterations per cycle $M$ may be comparable to the total number of
iterations needed to solve a particular problem for a prescribed
tolerance. At this point, using a method with higher $P$, and thus
higher $M$, would increase the number of iterations to converge, even
if the $\Gamma(\kappa)$ is nominally smaller.  With this limitation in
mind we outline a procedure to obtain optimal SRJ schemes, minimizing
the total number of iterations needed to reduce the residual by an
amount sufficient to reach convergence or, equivalently, to minimize
$|G_M(\kappa)|$. Note that the total number of iterations can be
chosen to be equal to $M$ without loss of generality, i.e. one cycle
of $M$ iterations is needed to reach convergence. To follow this
procedure one should find the optimal scheme for fixed values of $M$,
and then choose $M$ such that the maximum value of $|G_M(\kappa)|$ is
similar to the residual reduction needed to solve a particular
problem.  The first step, the minimization problem, is in general
difficult to solve, since fixing $M$ gives an enormous freedom in the
choice of the number of weights $P$, which can range from $1$ to $M$.
However, the numerical results of \cite{Yang2014} and
\cite{Adsuetal15}, seem to suggest that in general increasing the
number of weights $P$ will always lead to better convergence
rates. This leads us to conjecture that the optimal SRJ scheme, for
fixed $M$, is the one with $P=M$, i.e. all weights are different and
each weight is used once per cycle, $q_i=1,\, (i=1,\ldots,M)$. In
terms of the total amplification factor $G_M(\kappa)$, it is quite
reasonable to think that if one maximizes the number of different
roots by choosing $P=M$, the resulting function is, on average, closer
to zero than in methods with smaller number of roots, $P<M$, and one
might therefore expect smaller maxima for the optimal set of
coefficients.  One of the aims of this work is to compute the optimal
coefficients for this particular case and demonstrate that $P=M$ is
indeed the optimal case. 

Another goal of this paper is to show the performance of
optimal SRJ methods compared with optimal SOR algorithms applied to a
number of different discretizations of the Laplacian operator in
two-dimensional (2D) and three-dimensional (3D) applications
(Sect.\,\ref{sec:numexamples}). We will show that optimal SRJ methods
applied to high-order discretizations of the Laplacian, which yield 
iteration matrices that cannot be consistently ordered, perform very
similarly to optimal SOR schemes (when an optimal SOR weight can be
computed). We will further discuss that the trivial parallelization of the
SRJ methods outbalances the slightly better scalar performance of SOR
in some cases (Sect.\,\ref{sec:9-17discret}). Also, we will show that
the optimal weight of the SOR method can be suitably approximated by
functions related to the geometric mean of the set of weights obtained
for optimal SRJ schemes. This is of particular relevance when the
iteration matrix is non-consistently ordered and hence, the analytic
calculation of the optimal SOR weight is extremely intricate.   

\section {Optimal $P=M$ SRJ scheme}
\label{sec:OptimalCheb}

Let us consider a SRJ method with $P=M$ and hence $q_n=1,\,
(n=1,\ldots,M)$. For this particular choice, the amplification factor
$G_M(\kappa)$ is a polynomial of degree $M$ in $\kappa$ with $M$
different roots.  In this case, the set of weights $\omega _{n}$ that
minimizes the value of the maximum of $|G_M(\kappa)|$, given by
Eq.~(\ref{eq:gm}), in the interval $\kappa \in [\kmin,\kmax]$,
$0<\kmin\le \kmax$, can be determined by the following $M$ conditions:
\begin{equation}
G_M(0) = 1 \quad ; \quad
G_M(\kappa_n) = - G_M(\kappa_{n+1}), \quad n=0,\ldots,M-1, \\
\label{eq:conditions}
\end{equation}
where $\kappa_0 = \kmin$, $\kappa_M = \kmax$, and 
$\kappa_n, \text{ }  n=1,\ldots,M-1$ are the relative extrema of the function 
$G_M(\kappa)$. To simplify further we rescale $\kappa$ as follows:
\begin{equation}
\tilde{\kappa} = 2 \frac{\kappa-\kmin}{\kmax-\kmin} - 1. 
\end{equation}
As a function of $\tilde{\kappa}$ the amplification factor is
$\tilde{G}_M(\tilde{\kappa}) = G_M(\kappa(\tilde{\kappa}))$.
In the resulting interval, $\tilde{\kappa}\in[-1,1]$, there is a unique
polynomial of degree $M$ such that the absolute value of $\tilde{G}_M(\tilde{\kappa})$
at the extrema $\tilde{\kappa}_i$ is the same (fulfilling the last $M-1$
Eqs.~(\ref{eq:conditions})) and such that $\tilde{G}_M(\tilde{\kappa}(0))=1$.
This polynomial is proportional to the Chebyshev polynomial
of first kind of degree $M$, $T_M(\kappa)$, which can be defined through the identity
$T_M (\cos \theta) = \cos (M \,\theta)$. This polynomial satisfies that
\begin{gather}
|T_M(-1)| = |T_M(\tilde{\kappa}_n)| = |T_M(+1) |= 1, \quad n=1,\ldots,M-1,
\end{gather}
with $\tilde{\kappa}_i$ being the local extrema of
$T_M(\tilde{\kappa})$ in $[-1,1]$.  The constant of proportionality
can be determined requiring (Eq.\,\ref{eq:conditions}) $G_M(0)=1$, and
the amplification factor reads in this case
\begin{gather}
\tilde{G}_M(\tilde{\kappa}) =  \frac{T_M(\tilde{\kappa})}{T_M(\tilde{\kappa}(0))}  \quad; \quad
\tilde{\kappa}(0) = -\frac{(1+\kmin/\kmax)}{(1-\kmin/\kmax)} < -1.
\label{eq:Gequ}
\end{gather}
This result is equivalent to  Markoff's theorem\footnote{For an
  accesible proof of the original theorem \cite{Markoff:1916} , see Young's textbook \cite{Young:1971}, Theorem 9-3.1.}.
Note that the value of $\tilde{\kappa}(0)$ does not depend on the actual values
of $\kmin$ and $\kmax$, but only on the ratio $\kmin/\kmax$.
The roots and local extrema of the polynomial $T_M(\tilde{\kappa})$ are located, respectively, at
\begin{gather}
\tilde{\omega}_{n}^{-1} =-\cos\left(\pi \frac{2n-1}{2M}\right), \;\; n=1,\ldots,M, \\
\tilde{\kappa}_{n} = \cos\left(\pi\frac{n}{M}\right), \;\; n=1,\ldots,M-1, 
\end{gather}
which coincide with those of $\tilde{G}_M(\tilde{\kappa})$. Therefore, the set of weights
\begin{gather}
\omega_{n} = 2 \left[ \kmax+\kmin  
-\left(\kmax-\kmin\right) \cos\left(\pi\frac{2n-1}{2M}\right) \right ]^{-1},\; 
n=1,\ldots,M,
\label{eq:omegan}
\end{gather}
corresponds to the optimal SRJ method for $P=M$.

We have found with the simple analysis of this section that the optimal SRJ scheme when $P=M$ is fixed
turns out to be closely related to a Chebyshev iteration or Chebyshev semi-iteration for the solution of systems of
linear equations (see, for instance, \cite{Gutknecht:2002} for a review). This is especially easy to realize if we
consider the original formulation of this kind of methods, which appeared in the literature as special implementations
of the non-stationary or semi-iterative Richardson's method (RM, hereafter; see, e.g., \cite{Young:1953,Frank:1960} for
generic systems of linear equations, or \cite{Shortley:1953} for the application to boundary-value problems).
\cite{Yang2014} argued that, for a uniform grid, Eq.\,\ref{eq:SRJ} is identical to that of the RM
\cite{Richardson11}. There is, nevertheless, a minor difference between Eq.\,\ref{eq:SRJ} of the SRJ method and the RM
as it has been traditionally written \cite{Young:1954b}, that using our notation would be $u^{n+1} = u^n +
\hat{\omega}_{n} ( b - A u^n )$, which gives the obvious relation $\hat{\omega}_{n} =\omega_{n} d^{-1}$, in the case in
which all elements in $D$ are the same and equal to $d$. We note that this difference disappears in more modern
formulations of the RM (e.g., \cite{Opfer:1984}), in which the RM is also written as a fix point iteration of the form
$u^{n+1}=Tu^n+c$, with $T=I-M^{-1}A$, $c=M^{-1}b$ and $M$ any non-singular matrix. Differently from the RM in its
definition by Young \cite{Young:1954b}, our method in the case $M=1$ would fall in the category of stationary
Generalized Richardson's (GRF) methods according to the textbook of Young \citep[][chap.\,3]{Young:1971}. GRF methods
are defined by the updating formula
\begin{gather}
u^{n+1} = u^n + P(Au^n-b)
\end{gather}
where $P$ is any non singular matrix (in our case, $P=-\omega_{n}D^{-1}$).  In the original work of Richardson
\cite{Richardson11}, all the values of $\hat{\omega}_{n}$ where set either equal or evenly distributed in
$[a,b]$, where $a$ and $b$ are, respectively, lower and upper bounds to the minimum and maximum eigenvalues,
$\lambda_i$ of the matrix $A$ (optimally, $a=\min{(\lambda_i)}$, $b=\max{(\lambda_i})$). If a single weight is used
throughout the iteration procedure, a convenient choice is $\hat{\omega}=2/(b+a)$.\footnote{In the case of SRJ
schemes with $P=M$, it is easy to demonstrate (see \ref{sec:properties}) that the harmonic mean of the weights 
$\omega_n$ very approximately equals the value of the inverse weight of the stationary RM 
($2d^{-1}/(\kmax+\kmin)\simeq 2/(b+a)$).}

Yang \& Mittal \cite{Yang2014} state that the SRJ approach to maximizing
convergence is fundamentally different from that of the stationary RM. They argue that the RM aims to reduce
$\Gamma(\kappa)$ uniformly over the range $[\kmin,\kmax]$ by generating equally spaced nodes of $\Gamma$ in this
interval, while SRJ methods set a min-max problem whose goal is to minimize $|\Gamma|_{\rm max}$.\footnote{We note that
  this argument does not hold in the implementation of the non-stationary RM method made by Young \cite{Young:1953},
  since in this case one also attempts to minimize $|\Gamma|_{\rm max}$.} As a result, SRJ methods require computing a
set of weights yielding two differences with respect to the non-stationary RM in its original formulation \cite{Yang2014}:

\begin{enumerate}
\item the nodes in the SRJ method are not evenly
  distributed in the range $[\kmin,\kmax]$;

\item optimal SRJ schemes naturally have many repetitions of the same
  relaxation factor whereas RM generated distinct values of $\hat{\omega}_{n}$
  in each iteration of a cycle.
  \label{point2}
  \end{enumerate}
From these two main differences, \cite{Yang2014} conclude that while optimal SRJ schemes actually gain in convergence
rate over Jacobi method as grids get larger, the convergence rate gain for Richardson's procedure (in its original
formulation) never produces acceleration factors larger than 5 with respect to the Jacobi method. This result was
supported by Young in his Ph.D. thesis~\cite[][p.\,4]{Young:1950}, but on the basis of employing orderings of the
weights which did pile-up roundoff errors, preventing a faster method convergence (see point 2 below).

The difference outlined in point 1 above is non existent for GRF methods, where the eigenvalues of $A$ are not necessarily
evenly distributed in the spectral range of matrix A (i.e., in the interval $[a,b]$).  We note that Young
\cite{Young:1953} attempted to chose the $\hat{\omega}_{n}$ parameters of the RM to be the reciprocals of the roots of
the corresponding Chebyshev polynomials in $[a,b]$, which resulted in a method that is {\em almost the same} as
ours, but with two differences:

{\em First}, we do not need to compute the maximum and minimum eigenvalues of the matrix $A$; instead, we compute $\kmax$ and
  $\kmin$, which are related to the maximum and minimum frequencies that can be developed on the grid of choice
  employing an straightforward von Neumann analysis. Indeed, this procedure to estimate the maximum and minimum
    frequencies for the elliptic operators (e.g., the Laplacian) in the continuum limit allows applying it to matrices
    that are not necessarily consistently ordered, like, e.g., the ones resulting from the 9-point discretization of
    the Laplacian \cite{Adams:1988}. In Sect.\,\ref{sec:9-17discret} we  show how our method can be straightforwardly
    prescribed in this case and other more involved (high-order) discretizations of the Laplacian.

{\it Second}, in Young's method \cite{Young:1953} the two-term recurrence relation given by Eq.\,\ref{eq:SRJ} turned out to be
    unstable. Young found that the reason for the instability was the build up of roundoff errors in the evaluation of
    the amplification factor (Eq.\,\ref{eq:gm}), which resulted as a consequence of the fact that many of the values of
    $\omega_{n}$ can be much larger than one. Somewhat unsuccessfully, Young \cite{Young:1953} tried different orderings
    of the sequence of weights $\omega_{n}$, and concluded that, though they ameliorated the problem for small values of
    $M$, did not cure it when $M$ was sufficiently large. Later, Young \cite{Young:1954b,Young:1956} examines a number
      of orderings and concluded that some gave better results than others. However, the key problem of existence of
      orderings for which RM defines a stable numerical algorithm amenable to a practical implementation was not shown
      until the work of Anderssen \& Golub \cite{Anderssen:1972}. These authors showed that employing the ordering
      developed by Lebedev \& Finogenov \cite{Lebedev:1971} for the iteration parameters in the Chebyshev cyclic
      iteration method, the RM devised by Young \cite{Young:1953} was stable against the pile-up of round-off
      errors. However, Anderssen \& Golub \cite{Anderssen:1972} left open the question of whether other orderings are
      possible. In our case, numerical stability is brought about by the ordering of the weights in the iteration
    procedure. This ordering is directly inherited from the SRJ schemes of \cite{Yang2014}, and notably differs from the
    prescriptions given for two- or three-term iteration relations in Chebyshev semi-iterations \cite{Gutknecht:2002}
    and from those suggested by \cite{Young:1953}. Indeed, the ordering we use differs from that of
      \cite{Lebedev:1971,Nikolaev:1972,Lebedev:2002} (see \ref{sec:ordering}). Thus, though we do not have a theoretical
      proof for it, we empirically confirm that other alternative orderings work.

Taking advantage of the analysis made by \cite{Young:1953}, we
  point out that the average rate of convergence of the method in a
  cycle of $M$ iterations is
\begin{equation}
R_M = \frac{1}{M} \log{|T_M(\tilde\kappa(0))|},
\label{eq:Rp}
\end{equation} 
and it is trivial to prove that for $\kappa\in[\kmin,\kmax]$
\begin{equation}
G_M(\kappa) \leq \left | \frac{1}{T_M(\tilde{\kappa}(0))} \right | < 1,
\end{equation}
providing a simple way to compute an upper bound for the amplification
factor for the optimal scheme. This condition also guarantees the
convergence of the optimal SRJ method. Therefore, if we aim to reduce
the initial residual of the method by a factor $\sigma$, we have to
select a sufficiently large $M$ such that
\begin{equation}
\sigma \geq |T_M(\tilde\kappa(0))|^{-1}
\label{eq:Pmax}
\end{equation}

It only remains to demonstrate that the optimal SRJ scheme with $P=M$ is also the optimal 
SRJ scheme for any $P\le M$. Markoff's theorem states that for any polynomial $Q(x)$ of degree smaller or equal to $M$, such that $\exists x_0\in\mathbb{R}, x_0<-1$, with $Q(x_0)=1$, and $Q(x) \neq  T_M(x)/T_M(x_0)$, then
\begin{align}
\max{|Q(x)|} > \max{\left|\frac{T_{M}(x)}{T_M(x_0)}\right|} \quad \forall x \in [-1,1].
\end{align}
This theorem implies that any other polynomial of order $P\le M$, different from Eq.~(\ref{eq:Gequ}), is a poorer choice
as amplification factor.  The first implication is that $G_M(\tilde{\kappa}(0))<G_{M-1}(\tilde{\kappa}(0))$, i.e.,
increasing $M$ decreases monotonically the amplification factor $G_M(\kappa)$.  As a consequence, the per iteration
amplification factor $\Gamma_M (\kappa)$ also decreases by increasing $M$.  The second consequence is that the case
$P<M$ results in an amplification factor with larger extrema than the optimal $P=M$ case, and hence proves that our
numerical scheme leads to the optimal set of weights for any SRJ method with $M$ steps.  This confirms our intuition
that adding additional roots to the polynomial would decrease the value of its maxima, resulting in faster numerical
methods.  Though the SRJ algorithm with $P=M$ we have presented here turns out to be nearly equivalent to the
non-stationary RM of Young \cite{Young:1953}, in order to single it out as the optimum among the SRJ schemes, we will
refer to it as the {\em Chebyshev-Jacobi} method (CJM) henceforth.



\section{Numerical examples}
\label{sec:numexamples}

\subsection{Laplace equation}

In order to assess the performance of the new optimal set of schemes devised, we resort to the same prototype numerical example
  considered in \cite{Yang2014}, namely, the solution of the Laplace equation with homogeneous Neumann boundary conditions in
two spatial dimensions, in Cartesian coordinates and over a domain with unitary size:
\begin{equation}
\begin{cases}
\displaystyle\frac{\partial^2}{\partial x^2} u(x,y) + \frac{\partial^2}{\partial y^2} u(x,y) = 0, & (x,y) \in (0,1) \times (0,1) \\*[0.3cm]
\displaystyle\left.\frac{\partial}{\partial x} u(x,y) \right|_{x=0} =  \left.\frac{\partial}{\partial x} u(x,y) \right|_{x=1}  = 0, & y \in (0,1) \\*[0.4cm]
\displaystyle\left.\frac{\partial}{\partial y} u(x,y) \right|_{y=0} = \left.\frac{\partial}{\partial y} u(x,y) \right|_{y=1}  = 0, & x \in (0,1). \\
\end{cases}
\label{eq:Laplace}
\end{equation}
We consider a spatial discretization of the Laplacian operator employing a second-order, 5-point formula
\begin{equation}
\Delta u_{ij} = u_{i-1,j}+u_{i+1,j}+u_{i,j-1}+u_{i,j+1}-4u_{ij} = 0.
\label{eq:5-points}
\end{equation}
To compare the performance of different numerical schemes we 
monitor the evolution of the difference between two consecutive approximations
of the solution for the model problem specified in Eq.~(\ref{eq:Laplace}),
\begin{equation}
||r^{n}||_{\infty} = \max_{ij} {|u^{n}_{ij} - u^{n-1}_{ij}|},
\label{eq:residual}
\end{equation}
where $u^n_{ij}$ is the numerical approximation computed after $n$ iterations at the grid point $(x_i,y_j)$.
\begin{figure}[h]
\centering
\includegraphics[width=0.49\textwidth]{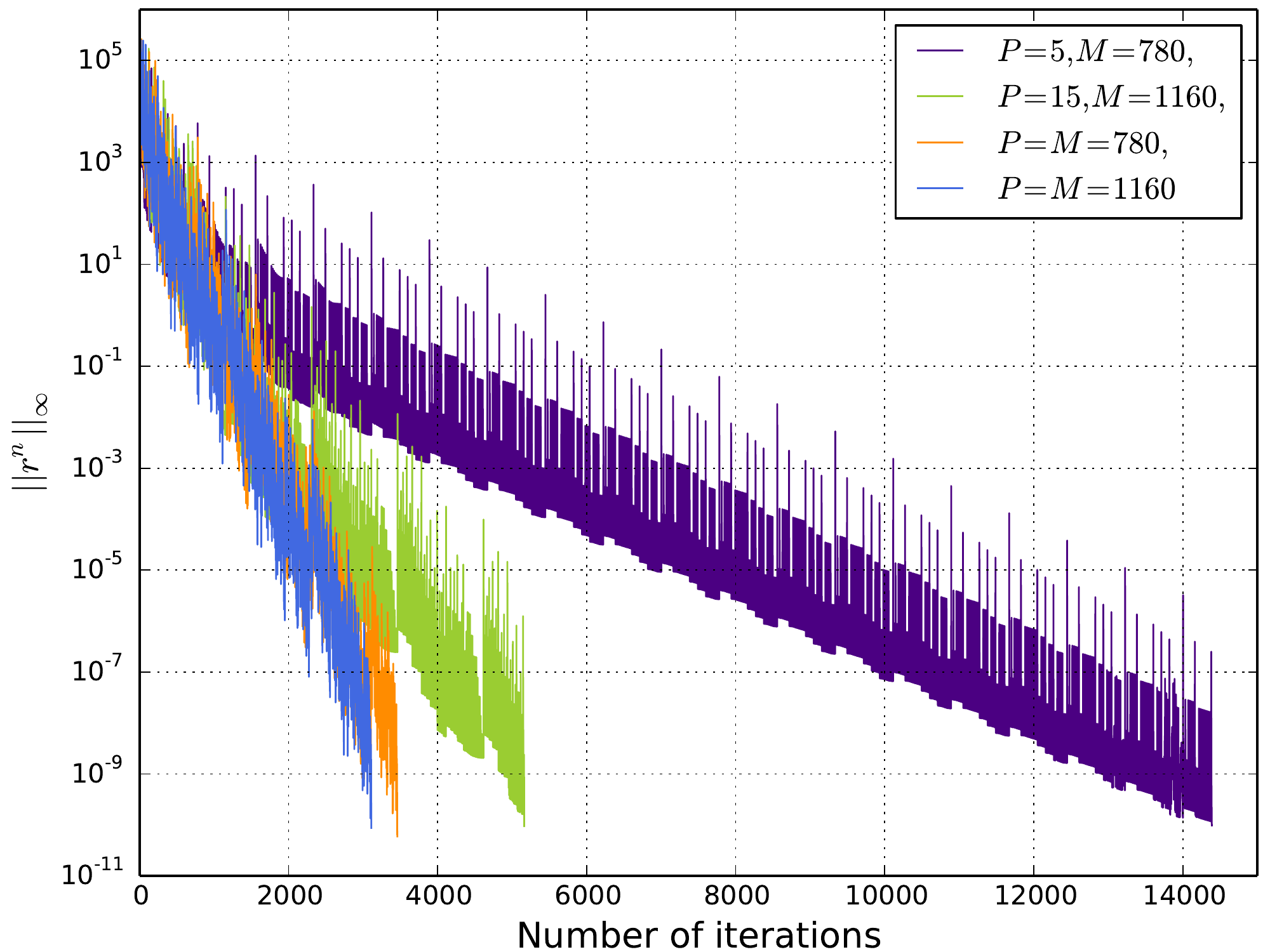}
\includegraphics[width=0.49\textwidth]{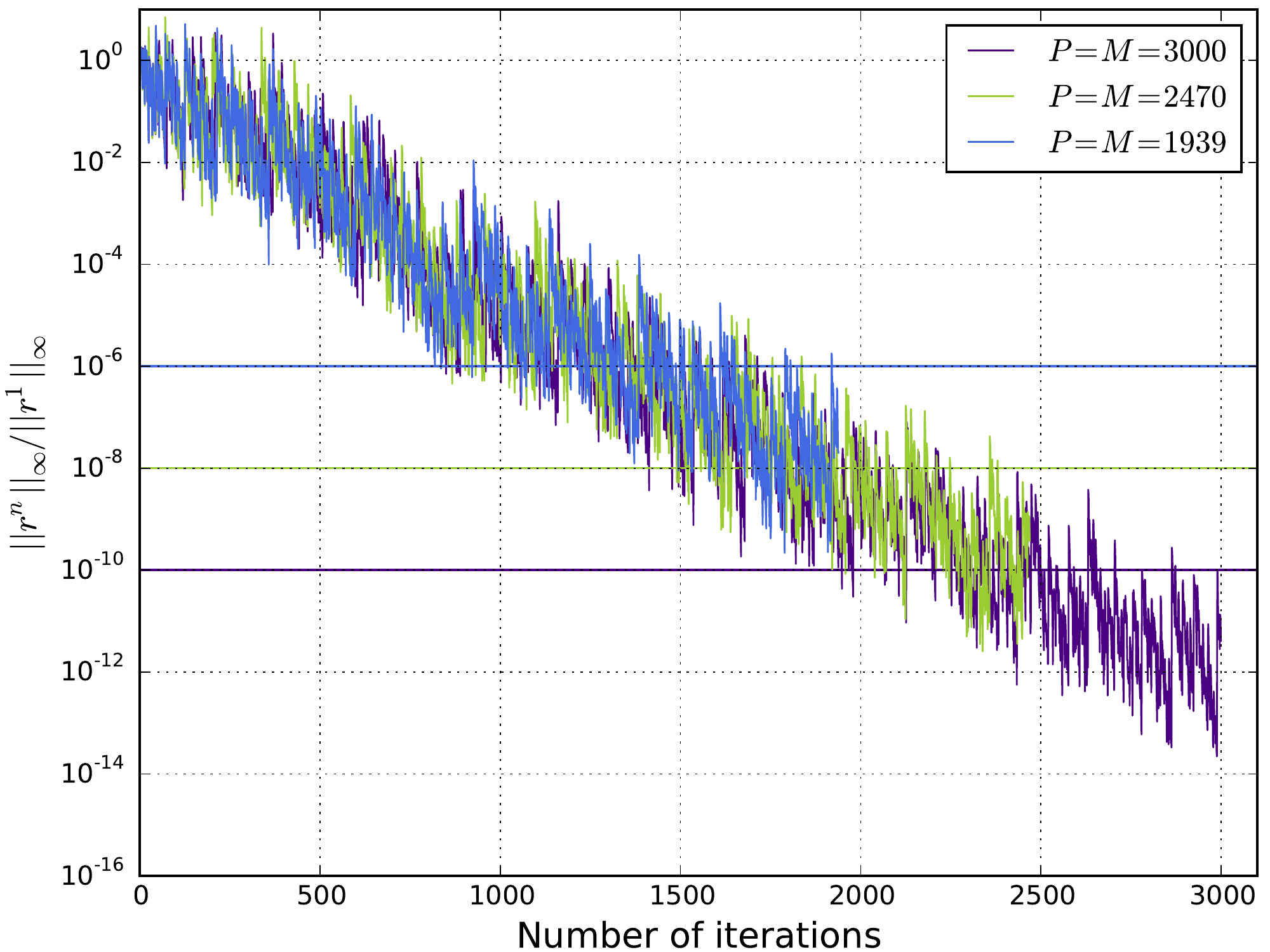} 
\caption{Left: Evolution of the residual $||r^n||_\infty$, defined in
  Eq.\,\ref{eq:residual}, as a function of the number of iterations
  for the problem set in Eq.\,\ref{eq:Laplace} and a Cartesian grid of
  $256\times 256$ uniform zones. The different color lines correspond
  to different schemes. Violet line: SRJ method with $P=5$ and
  $M=780$. Orange line: CJM with $P=M=780$. Green line: SRJ method
  with $P=15$ and $M=1160$. Blue line: CJM with $P=M=1160$. We can
  observe that the reduction of the residual is faster in the new
  Chebyshev-Jacobi schemes than in the corresponding SRJ schemes with
  the same value of $M$. Right: We show three examples where we
  computed the optimal value of the $M$ for reaching the desired
  residual in one cycle. The cases $P=1939$, 2470 and 3000 correspond
  to schemes that (theoretically) should reduce the initial residual
  by factors $\simeq 10^6$, $10^8$ and $10^{10}$.}
\label{fig:fig05}
\end{figure}

\begin{figure}[htb]
\centering
\includegraphics[width=0.49\textwidth]{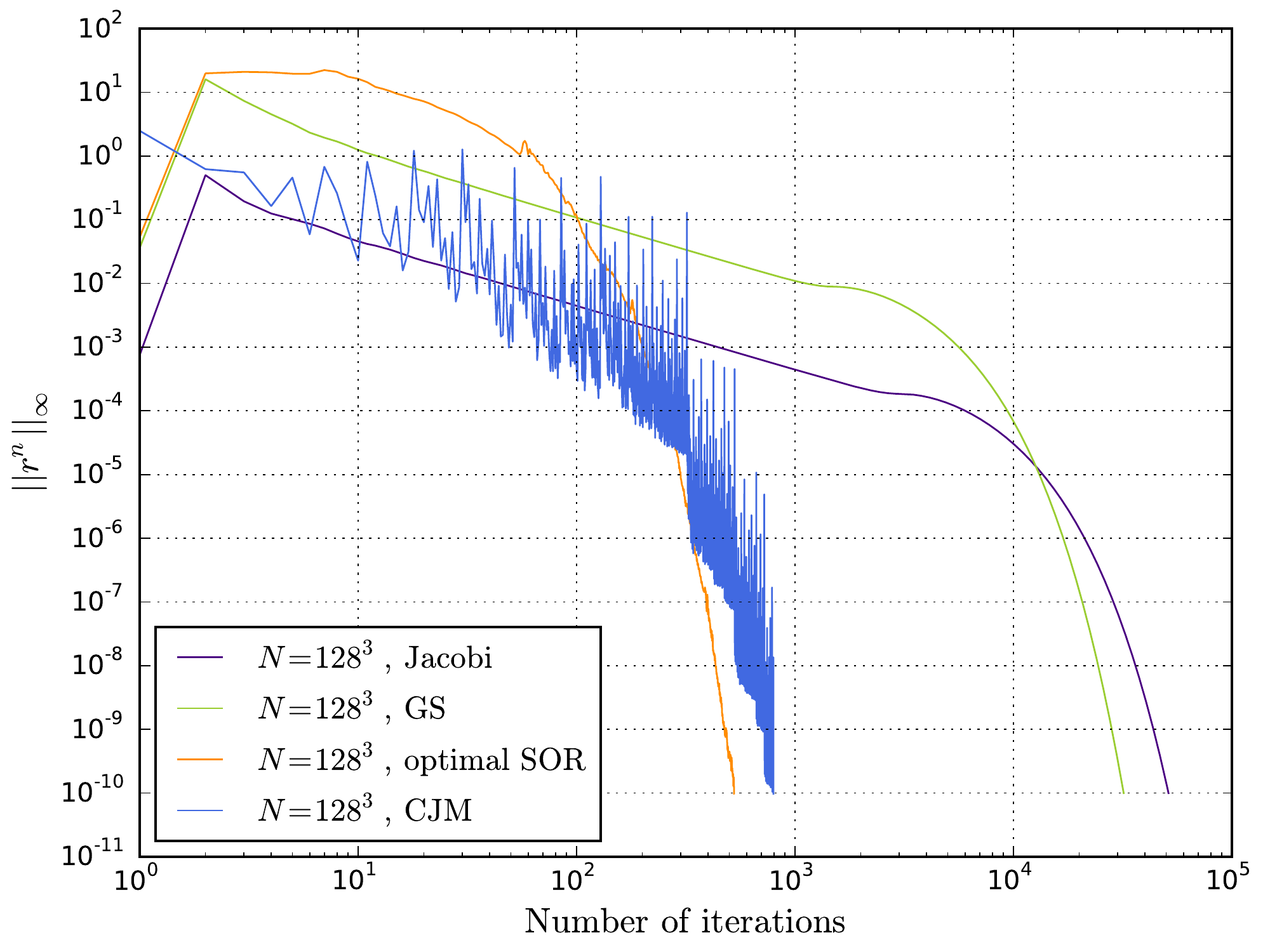} 
\includegraphics[width=0.49\textwidth]{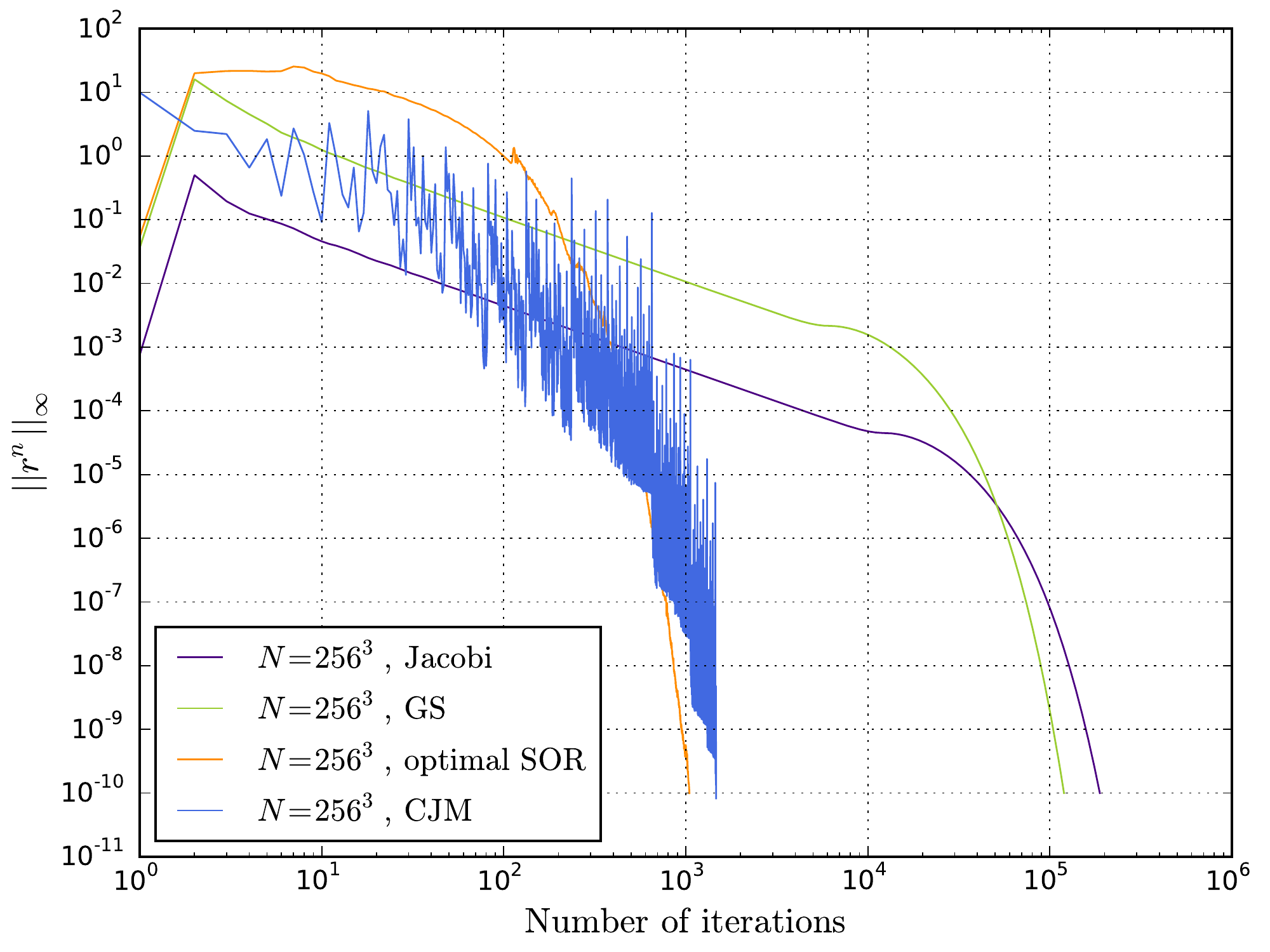}
\caption{The evolution of the residual for the solution of the Poisson equation~\eqref{eq::poisson3D} in 3D,
    with $N=128$ (left panel) and $N=256$ (right panel) for different iterative methods.} 
\label{fig:poisson3D}
\end{figure}

\begin{figure}[htb]
\centering
\includegraphics[width=0.75\textwidth]{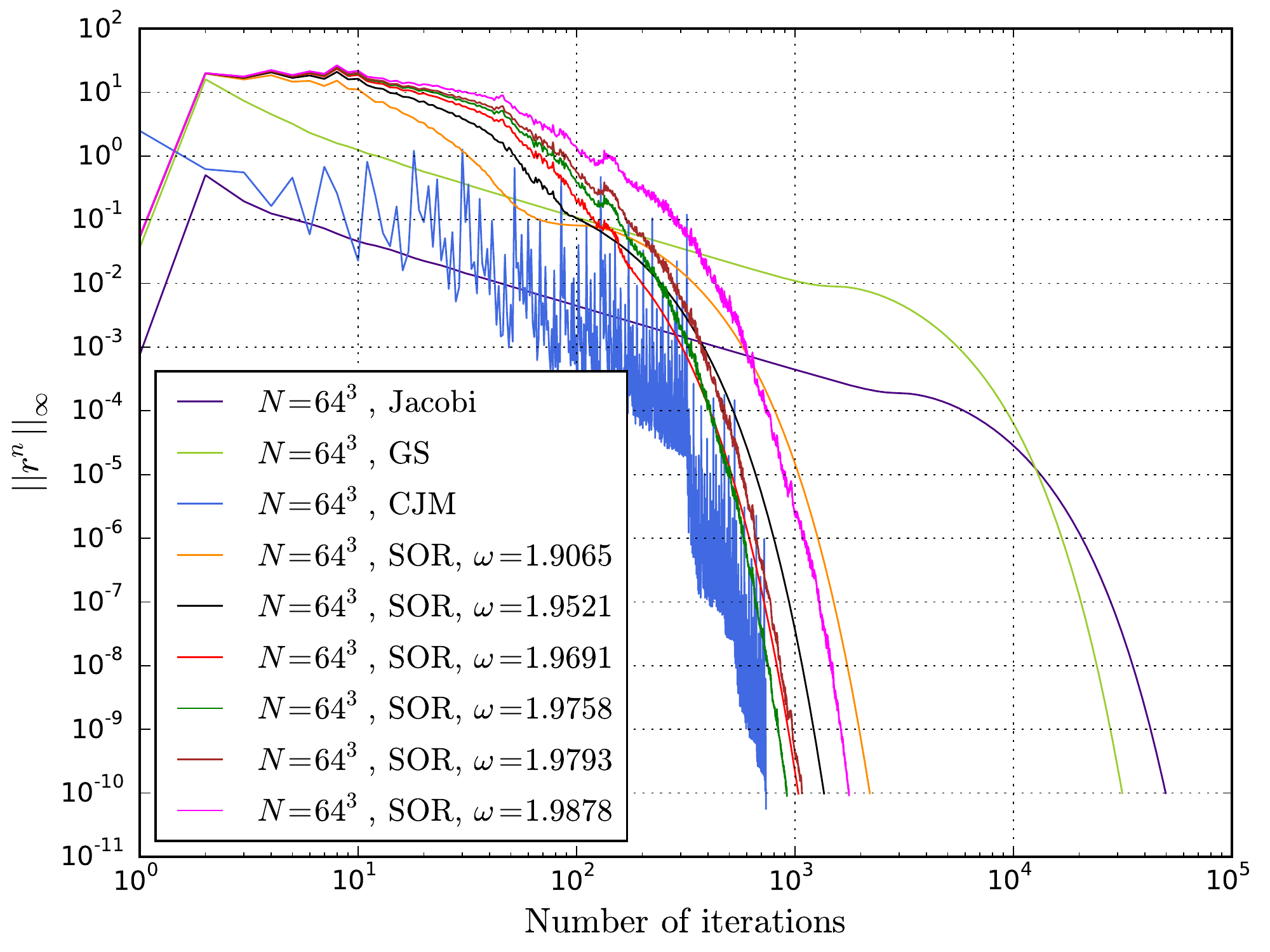} 
\caption{The evolution of the residual for the solution of the Poisson equation~\eqref{eq::poisson3D} in 3D
    using octant symmetry, with $N=64$ for different iterative methods and different relaxation factors $\omega$ in SOR.}
\label{fig:poisson3D_octant}
\end{figure}

In Fig.\,\ref{fig:fig05} (left), we compare the evolution of the residual as a function of the number of iterations for
several SRJ schemes, as well as for the new schemes developed here. The violet line corresponds to the best SRJ scheme
presented in \cite{Yang2014} for the solution of the problem set above and a spatial grid of $N_x \times N_y = 256
\times 256$ uniform zones, i.e. the SRJ scheme with $P=5$ and $M=780$. Comparing with the new CJM for $P=M=780$ (orange
line in Fig.\,\ref{fig:fig05} left), it is evident that the new scheme
reduces the number of iterations to
reach the prescribed tolerance ($||r^n||_\infty \le 10^{-10}$ in this
example) by about a factor of 5. We also include in Fig.\,\ref{fig:fig05}
(left; green line) the residual evolution corresponding to the best SRJ optimal algorithm developed by \cite{Adsuetal15}
for the proposed resolution, namely, the scheme with $P=15$ levels and $M=1160$. It is obvious that even the CJM with
$P=M=780$ reduces the residual faster than the $P=15$ SRJ scheme. However, since the $P=15$ SRJ scheme requires a larger
value of $M$ than in the case of $P=5$, for a fair comparison, we also include in Fig.\,\ref{fig:fig05} (left; blue
line) the CJM with $P=M=1160$. The latter is the best performing scheme, thought the difference between the two new CJM
with different values of $P$ is very small (in Fig.\,\ref{fig:fig05} the blue and orange lines practically overlap).

A positive property of the new algorithm presented in Sect.\,\ref{sec:OptimalCheb} is its predictability, i.e., the
easiness to estimate the size of the $M$-cycle in order to reduce the tolerance by a prescribed amount
(Eq.\,\ref{eq:Pmax}). Indeed, it is not necessary to monitor the evolution of the residual in every iteration (as in
many other non-stationary methods akin to the Richardson's method -e.g., in the gradient method-), with the obvious
reduction in computational load per iteration that this implies. In Fig.\,\ref{fig:fig05} (right) we show that our
algorithm performs as expected, reducing the initial residual by factors of larger than $10^6$, $10^8$ and $10^{10}$ in
a single cycle consistent of $P=1939$, 2470 and 3000 iterations, respectively, since for the problem at hand we have
$\kmin=\sin{(\frac{\pi}{2 \times 256})}^2=3.76491\times 10^{-5}$, $\kmax=2$, and thus, $\tilde\kappa(0)=-1.00004$.

In this simple example the upper bound for the residual obtained from Eq.\,(\ref{eq:Pmax}) is very rough and clearly
overestimates the number of iterations to reduce the residual below the prescribed values. In more complex problems this
will not necessarily be the case as we will show in the following, more demanding example.

\subsection{Poisson equation in 3D}

Here we test the CJM and the
predictability of the residual evolution in a three-dimensional
elliptic equation with a source term.  
For this test, we use infrastructure provided by the Einstein
Toolkit~\cite{ET,Loffler2012}.
The actual calculation is finding the static field of a uniformly charged sphere of radius $R$ in 3D Cartesian coordinates subject to Dirichlet
boundary conditions, solving the Poisson equation:
\begin{equation}\label{eq::poisson3D}
\Delta \phi(x,y,z) = -4 \pi \rho, 
\end{equation}
where $\rho = \frac{3Q}{4 \pi R^3}$ and $Q$ is the charge of the sphere. We solve the elliptic
equation~\eqref{eq::poisson3D} with a standard second-order accurate 7-point stencil
\begin{equation}
\Delta u_{ijk} = u_{i-1,jk}+u_{i+1,jk}+u_{i,j-1,k}+u_{i,j+1,k}+u_{i,j,k-1}+u_{ij,k+1}-6u_{ijk} = 0.
\label{eq:7-points}
\end{equation}

We consider two different grid sizes with $N_x=N_y=N_z=N=128$ and
$N_x=N_y=N_z=N=256$ points and the following iterative methods:
Jacobi, Gauss-Seidel (SOR with $\omega=1$), SOR with the optimal relaxation factor
$\omega_{\rm opt}=2/(1+\sin{(\pi/N)})$, and CJM with the optimal
sequence of weights for a given resolution.  The results for the two
grid resolutions are shown in Fig.~\ref{fig:poisson3D}. Both SOR and
CJM (slightly less than twice the number of iterations of SOR) are
more than an order of magnitude faster than the Jacobi and
Gauss-Seidel methods.  While the CJM method is not as fast as SOR when
using the optimal relaxation factor $\omega_{\mathrm{opt}}$, we note
here two arguments that should favor the use of the CJM over SOR:
Firstly, Young's theory of relating $\omega_{\mathrm{opt}}$ to the
spectral radius of the Jacobi iteration matrix $\rho(J)$ via
$\omega_{\mathrm{opt}} = 2/(1+\sqrt{1-\rho(J)^2})$ only applies when
the original matrix of the linear system $Au=b$ is consistently
ordered. Secondly, the CJM method is trivially parallelized, while SOR
requires multicolor schemes for a successful parallelization, as we
will discuss below presenting results for 9-point and 17-point
Laplacians in 2D.
 
Next, we solve equation~\eqref{eq::poisson3D} subject to reflection
symmetry (homogeneous Neumann boundary conditions) at the $x=0, y=0,
z=0$ planes (so-called octant symmetry) with $N_x=N_y=N_z=N=64$
points, using the same iterative methods as before. For the CJM, we
choose the same sequence of weights as those we used for the
full 3D domain using $N=128$ points. Because of the boundary
conditions used to impose octant symmetry, the resulting matrix $A$ is
non-consistently ordered and hence there is no analytic expression to
calculate $\omega_{\mathrm{opt}}$ for SOR; in this case we test a
sequence of values of $\omega$ to empirically estimate the optimal
value for the given problem. The residuals of the different iterative methods are shown in
Fig.~\ref{fig:poisson3D_octant}. The CJM now performs better than SOR
for any $\omega$ we have tested. Furthermore, as seen in the plot, SOR
is very sensitive to the exact value of $\omega$ that is chosen, as is
well known.  The CJM method is free of this need to estimate and choose a sensitive
parameter.


\subsection{CJM for non-consistently ordered matrices: high-order discretization of the Laplacian operator in 2D with 9 and 17 points}
\label{sec:9-17discret}

As we have already mentioned, Young's theory of relating the optimal
SOR parameter to the spectral radius of the Jacobi iteration matrix
does not apply in the case of non-consistently ordered matrices. In
this section, we will investigate two of these cases, namely a 9-point
and 17-point discretization of the Laplacian in 2D.
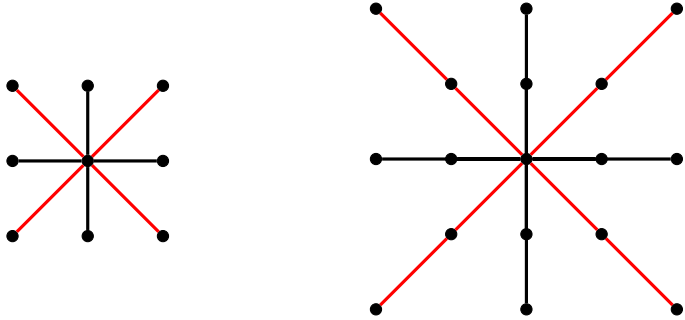
\begin{figure}
\centering
\begin{tikzpicture}[baseline=-58.5pt]
  \stencilpt{-1,0}{i-1}{};
\stencilpt{-1,1}{i-1--j+1}{};
\stencilpt{-1,-1}{i-1--j-1}{};
\stencilpt{1,-1}{i+1--j-1}{};
\stencilpt{1,1}{i+1--j+1}{};
  \stencilpt{ 0,0}{i}  {};
  \stencilpt{ 1,0}{i+1}{};
  \stencilpt{0,-1}{j-1}{};
  \stencilpt{0, 1}{j+1}{};
  \draw[very thick] (j-1) -- (i)
        (i)   -- (j+1)
        (i-1) -- (i)
        (i)   -- (i+1);
        \draw[red,very thick](i-1--j+1) -- (i)
        (i-1--j-1) -- (i)
        (i+1--j+1) -- (i)
        (i+1--j-1) -- (i);
\end{tikzpicture}
\qquad \qquad \qquad
\begin{tikzpicture}
  \stencilpt{-1,0}{i-1}{};
\stencilpt{-1,1}{i-1--j+1}{};
\stencilpt{-1,-1}{i-1--j-1}{};
\stencilpt{1,-1}{i+1--j-1}{};
\stencilpt{1,1}{i+1--j+1}{};
  \stencilpt{ 0,0}{i}  {};
  \stencilpt{ 1,0}{i+1}{};
  \stencilpt{0,-1}{j-1}{};
  \stencilpt{0, 1}{j+1}{};
\stencilpt{-2,0}{i-2}{};
\stencilpt{2,0}{i+2}{};
\stencilpt{0,-2}{j-2}{};
\stencilpt{0,2}{j+2}{};
\stencilpt{-2,2}{i-2--j+2}{};
\stencilpt{-2,-2}{i-2--j-2}{};
\stencilpt{2,-2}{i+2--j-2}{};
\stencilpt{2,2}{i+2--j+2}{};
  \draw[very thick] (j-1) -- (i)
        (i)   -- (j+1)
        (i-1) -- (i)
        (i)   -- (i+1)
        (j-2) -- (i)
        (i)   -- (j+2)
        (i-2) -- (i)
        (i)   -- (i+2);
       \draw[red,very thick] (i-1--j+1) -- (i)
        (i-1--j-1) -- (i)
        (i+1--j+1) -- (i)
        (i+1--j-1) -- (i)
        (i-2--j+2) -- (i-1--j+1)
        (i-2--j-2) -- (i-1--j-1)
        (i+2--j+2) -- (i+1--j+1)
        (i+2--j-2) -- (i+1--j-1);
\end{tikzpicture}
\caption{Schematic representation of the 9- and 17-point stencils. The
black and red lines correspond to the standard stencil $S_+$ and
rotated stencil $S_\times$, respectively. See main text for details.}
\label{fig::comp_molecules}
\end{figure}

One way of obtaining this type of discretizations is doing a convex
combination between the discretization of the Laplacian operatior
using the standard stencil, $S_+$, with its discretization in a
rotated stencil, $S_\times$ (see Fig.~\ref{fig::comp_molecules}):
\begin{gather}
\alpha S_+ + (1-\alpha) S_\times \, .
\label{eq:gendiscret}
\end{gather}  
Writing $\alpha$ as a rational number $a/b$, the resulting 9-points discretized Laplacian is
\begin{gather}
\Delta u_{ij}=\frac{1}{2b h^2} \bigg [ 2 a u_{i-1,j} + 2 a u_{i+1,j} + 2 a u_{i,j-1} + 2 a u_{i,j+1} \nonumber \\
+ (b-a) u_{i-1,j-1} + (b-a) u_{i+1,j+1} + (b-a) u_{i-1,j+1} + (b-a) u_{i+1,j-1} \nonumber \\
- 4(a+b) u_{i,j} \bigg ],
\label{eq:genLap} 
\end{gather}  
where, for simplicity, we assume that the grid spacing, $h$, is the same in the $x-$ and $y-$directions.
From this general form, we can recover the standard 5-points discretization simply taking $a=b=1$. In the same
way, we can recover the 9-points discretization of the Laplacian studied in \cite{Adams:1988} by imposing $a=2$ and
$b=3$:
\begin{gather}
\Delta u_{ij}=\frac{1}{6 h^2} \bigg [ 4 u_{i-1,j} + 4 u_{i+1,j} + 4u_{i,j-1} + 4 u_{i,j+1} \nonumber \\
+ u_{i-1,j-1} + u_{i+1,j+1} + u_{i-1,j+1} + u_{i+1,j-1} - 20 u_{i,j}
\bigg ]\, .
\label{eq:9-points}
\end{gather}  
From the von Neumann stability analysis of Eq.\,(\ref{eq:Laplace}), we obtain the following expression of the
amplification factor for the Laplacian discretization of Eq.\,(\ref{eq:genLap})
\begin{gather}
G = 1 - \omega \bigg [ \frac{2 a}{a+b} \sin^2{\frac{k_x \Delta x}{2}} + \frac{2 a}{a + b} \sin^2{\frac{k_y \Delta y}{2}}
\nonumber \\
+ \frac{b-a}{a+b} \big [ 1 - \cos{(k_x \Delta_x)} \cos{(k_y \Delta_y)}
\big ] \bigg ]
\label{eq:Gampl}
\end{gather}  
For $\alpha=a=b=1$, we recover the expression of the amplification factor shown in \cite{Yang2014,Adsuetal15}. It is
easy to check that when $a = 2$ and $b=3$, Eq.\,\eqref{eq:Gampl} reduces to
\begin{gather}
G= 1 - \frac{ \omega}{5}\bigg [4 \sin^2{\frac{k_x \Delta x}{2}} + 4
\sin^2{\frac{k_y \Delta y}{2}} + 1 - \cos{(k_x \Delta_x)}
\cos{(k_y \Delta_y)} \bigg ].
\end{gather}  
The factor multiplying $\omega$ in the previous expression is related to the weights of any SRJ scheme and singularly
with the CJM. As a function of the wave number $\kappa$, the minimum amplification factor results for $k_x = k_y = \pi$,
while the maximum amplification factor is attained for $k_x =\pi /\Delta_x$ and $k_y = \pi /\Delta_y$, with respective wave
numbers $\kmin$ and $\kmax$, whose expressions are
\begin{gather}
\kmin = \frac{4}{5} \sin^2{\frac{\pi}{2 N_x}} + \frac{4}{5} \sin^2{\frac{\pi}{2 N_y}} + \frac{1}{5} \bigg [ 1 - \cos{\frac{\pi}{Nx} } \cos{\frac{\pi}{N_y} } \bigg ] ,\\
\kmax = \frac{8}{5} \,.
\label{Eq:kmkM9p}
\end{gather}  
It can be shown that the 9-point discretization of the Laplacian
provides a fourth-order accurate method for the Poisson equation when
the source term is smooth \cite{leveque2007finite}. 


Next, we consider the case of a 17-point discretization of the Laplacian. From the general form of
Eq.~\eqref{eq:gendiscret}, again writing $\alpha=a/b$ one obtains
\begin{gather}
\Delta u_{ij} = \frac{1}{24 b h^2} \bigg [ \nonumber\\
- 2 a u_{i-2,j} + 32 a u_{i-1,j} + 32 a u_{i+1,j} - 2 a u_{i+2,j} \nonumber\\
- 2 a u_{i,j-2} + 32 a u_{i,j-1} + 32 a u_{i,j+1} - 2 a u_{i,j+2}  \nonumber\\
- (b-a) u_{i-2,j-2} + 16(b-a) u_{i-1,j-1} + 16 (b-a) u_{i+1,j+1} - (b-a) u_{i+2,j+2} \nonumber \\
- (b-a) u_{i-2,j+2} + 16(b-a) u_{i-1,j+1} + 16 (b-a) u_{i+1,j-1} - (b-a) u_{i+2,j-2} \nonumber\\
- 60 (a+b) u_{i,j} \bigg ].
\label{eq:genLap17} 
\end{gather}  
The standard 9-point discretization of the Laplacian is recovered for
$a=b=1$ in Eq.\,\eqref{eq:genLap17}. Performing the von Neumann stability
analysis for Eq.\,(\ref{eq:Laplace}), we obtain the following expression of the amplification factor for the Laplacian
discretization of Eq.\,(\ref{eq:genLap17})
\begin{gather}
G = 1-\omega\frac{1}{15(a+b)} \bigg [ 
-2a \big(\sin^2 (k_x \Delta_x)+ \sin^2 (k_y \Delta_y) \big) +  \nonumber\\  
32 a \bigg ( \sin^2 \bigg ( \frac{k_x \Delta_x}{2} \bigg ) + \sin^2 \bigg ( \frac{k_y \Delta_y}{2}\bigg ) \bigg) - \nonumber\\
(b-a)\big( [1 - \cos (2 k_x \Delta_x) \cos (2 k_y \Delta_y)] - 16 [1 - \cos (k_x \Delta_x) \cos (k_y \Delta_y)]\big )
\bigg ],
\end{gather}  
and,therefore, taking into account the minimum and maximum wave numbers as 
in the previous case, the extremal values of $\kappa$ are: 
\begin{gather}
\kmin =
\frac{1}{15(a+b)} \bigg[ -2a \bigg( \sin^2 \frac{\pi }{N_x} +  \sin^2 \frac{\pi }{N_y} \bigg) 
+ 32 a\bigg( \sin^2 \frac{\pi}{2 N_x}  + \sin^2 \frac{\pi}{2 N_y} \bigg) \nonumber\\
-(b-a) \bigg([1 - \cos \frac{2 \pi}{N_x} \cos \frac{2 \pi}{N_y}] - 16 [1 - \cos \frac{\pi}{N_y} \cos
\frac{\pi}{N_y}]\bigg)\bigg]
\label{eq:kmin17p},\\
\kmax = \frac{64 a}{15(a+b)} \,.
\label{eq:kmax17p}
\end{gather}  
Let us consider the particular case $a=1$ and $b=2$. For the Laplacian discretization \eqref{eq:genLap17}, we have
\begin{gather}
\Delta u_{ij}  = \frac{1}{48 h^2} \bigg [
- 2 u_{i-2,j} + 32 u_{i-1,j} + 32 u_{i+1,j} - 2 u_{i+2,j} \nonumber\\
- 2 u_{i,j-2} + 32 u_{i,j-1} + 32 u_{i,j+1} - 2 u_{i,j+2}  \nonumber\\
- u_{i-2,j-2} + 16 u_{i-1,j-1} + 16  u_{i+1,j+1} - u_{i+2,j+2} \nonumber\\ 
- u_{i-2,j+2} + 16 u_{i-1,j+1} + 16  u_{i+1,j-1} - u_{i+2,j-2} - 180 u_{i,j} \bigg ]
\label{eq:17-points}
\end{gather}  
and the expressions for $\kmin$ and $\kmax$ of Eqs.\,\eqref{eq:kmin17p} and \eqref{eq:kmax17p} reduce to
\begin{gather}
\kmin =
\frac{1}{45} \bigg[ -2\bigg(\sin^2 \frac{\pi }{N_x} + \sin^2 \frac{\pi }{N_y} \bigg) + 32\bigg( \sin^2 \frac{\pi}{2
  N_x}  + \sin^2 \frac{\pi}{2 N_y} \bigg)\nonumber\\
- [1 - \cos \frac{2 \pi}{N_x} \cos \frac{2 \pi}{N_y}] + 16 [1 - \cos \frac{\pi}{N_y} \cos \frac{\pi}{N_y}] \bigg]\\
\kmax = \frac{64}{45} 
\end{gather}  

\begin{figure}[h]
\centering
\includegraphics[width=0.49\textwidth]{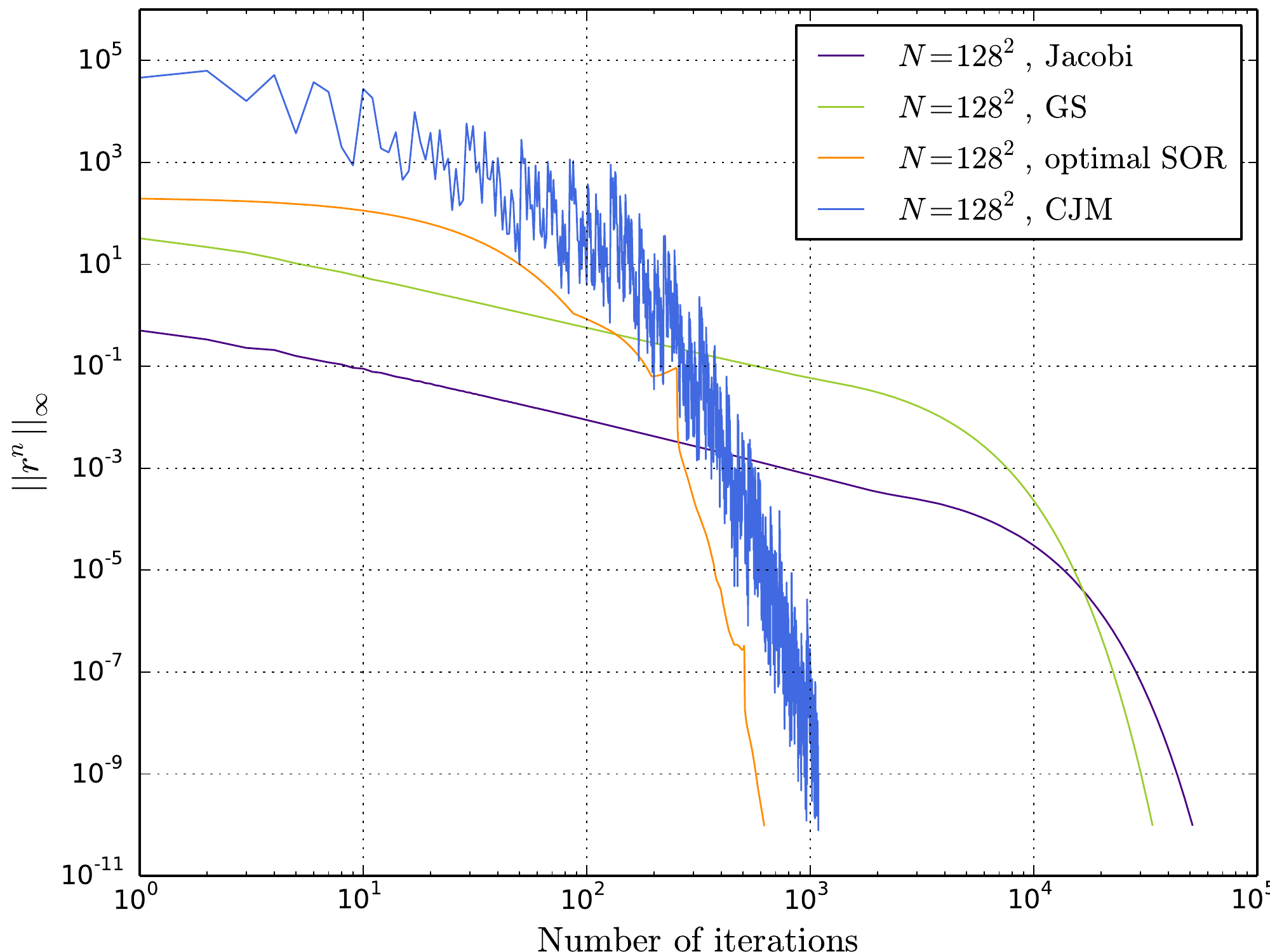} 
\includegraphics[width=0.49\textwidth]{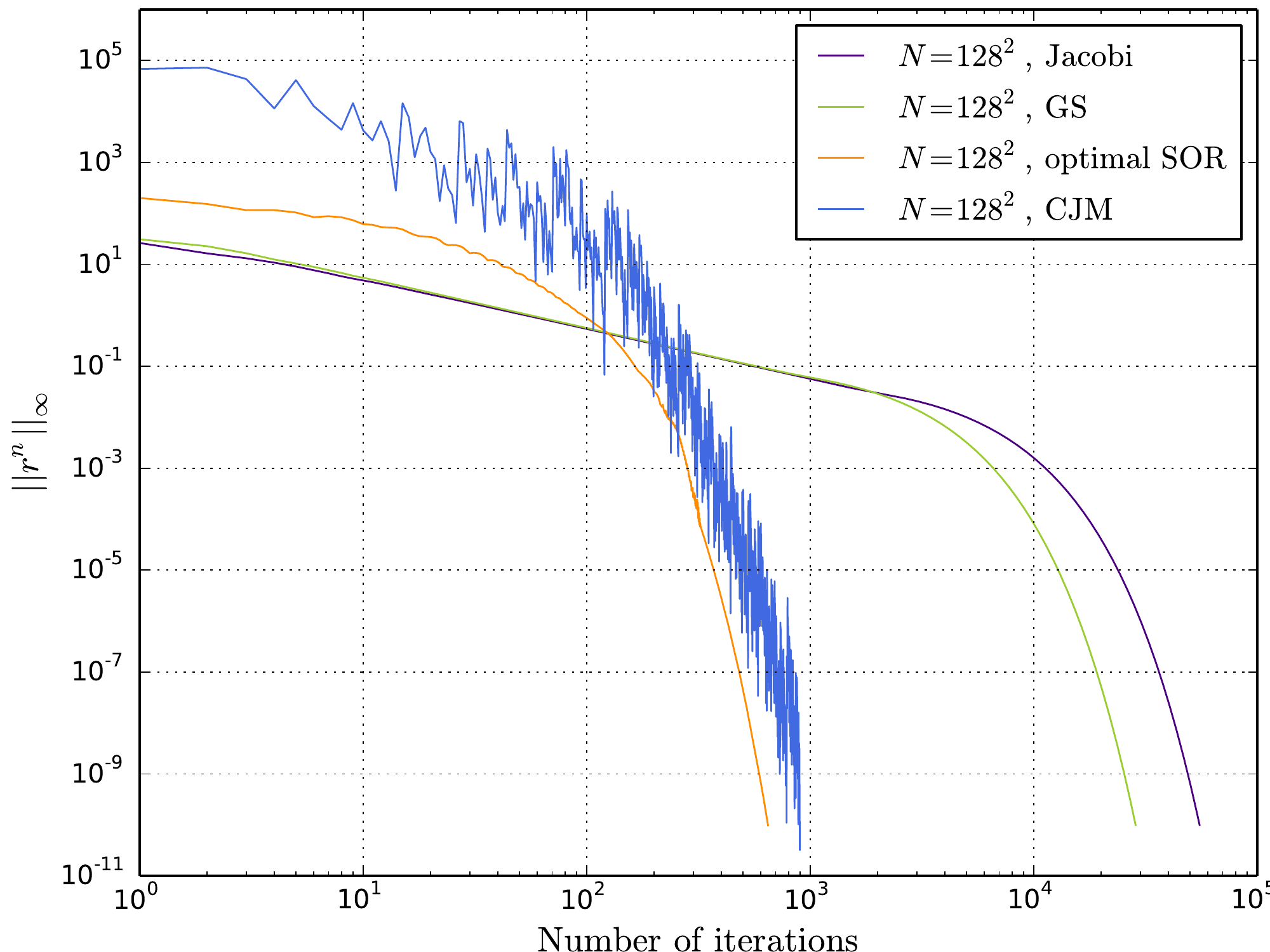} \\
\includegraphics[width=0.49\textwidth]{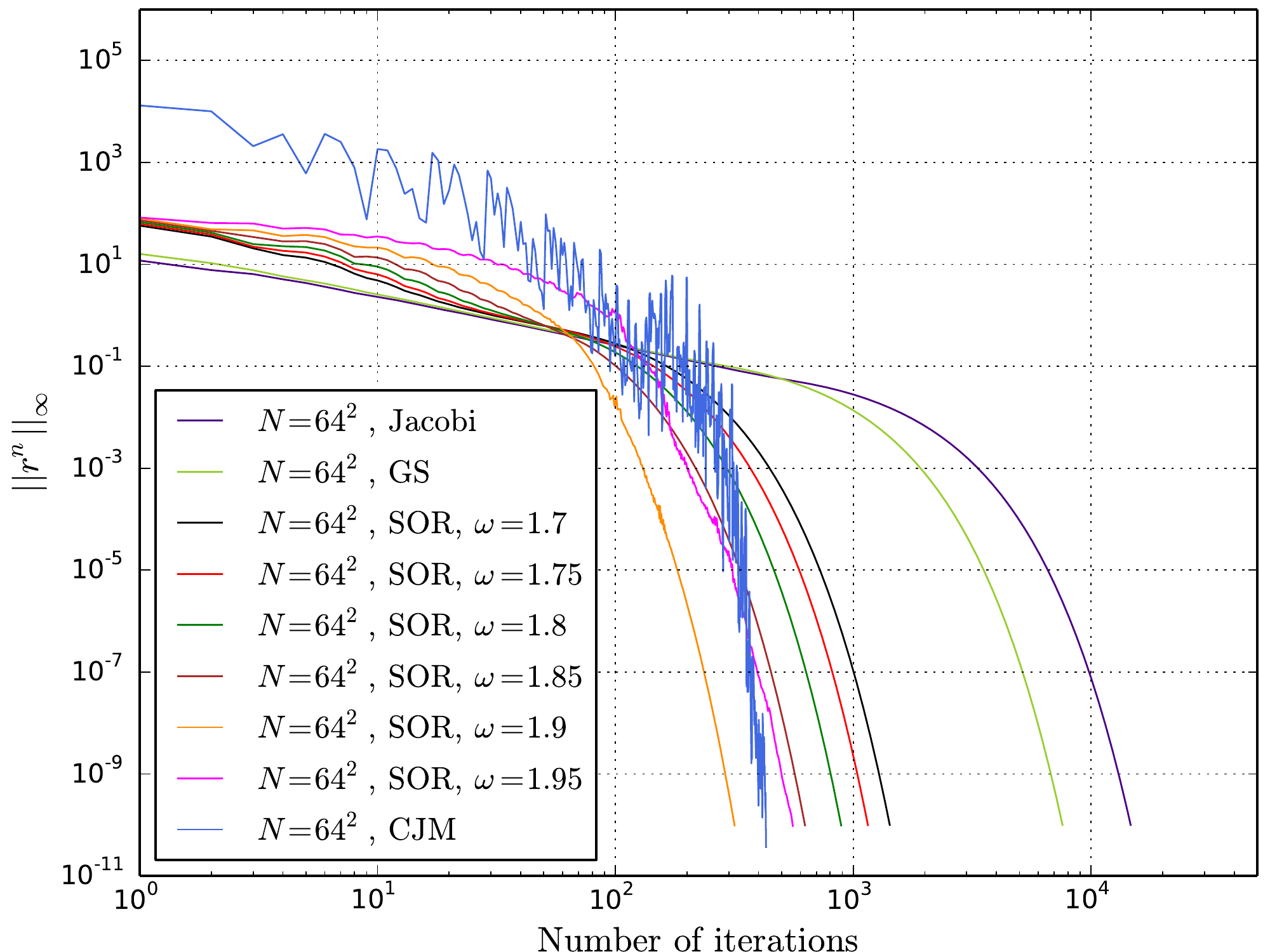} 
\includegraphics[width=0.49\textwidth]{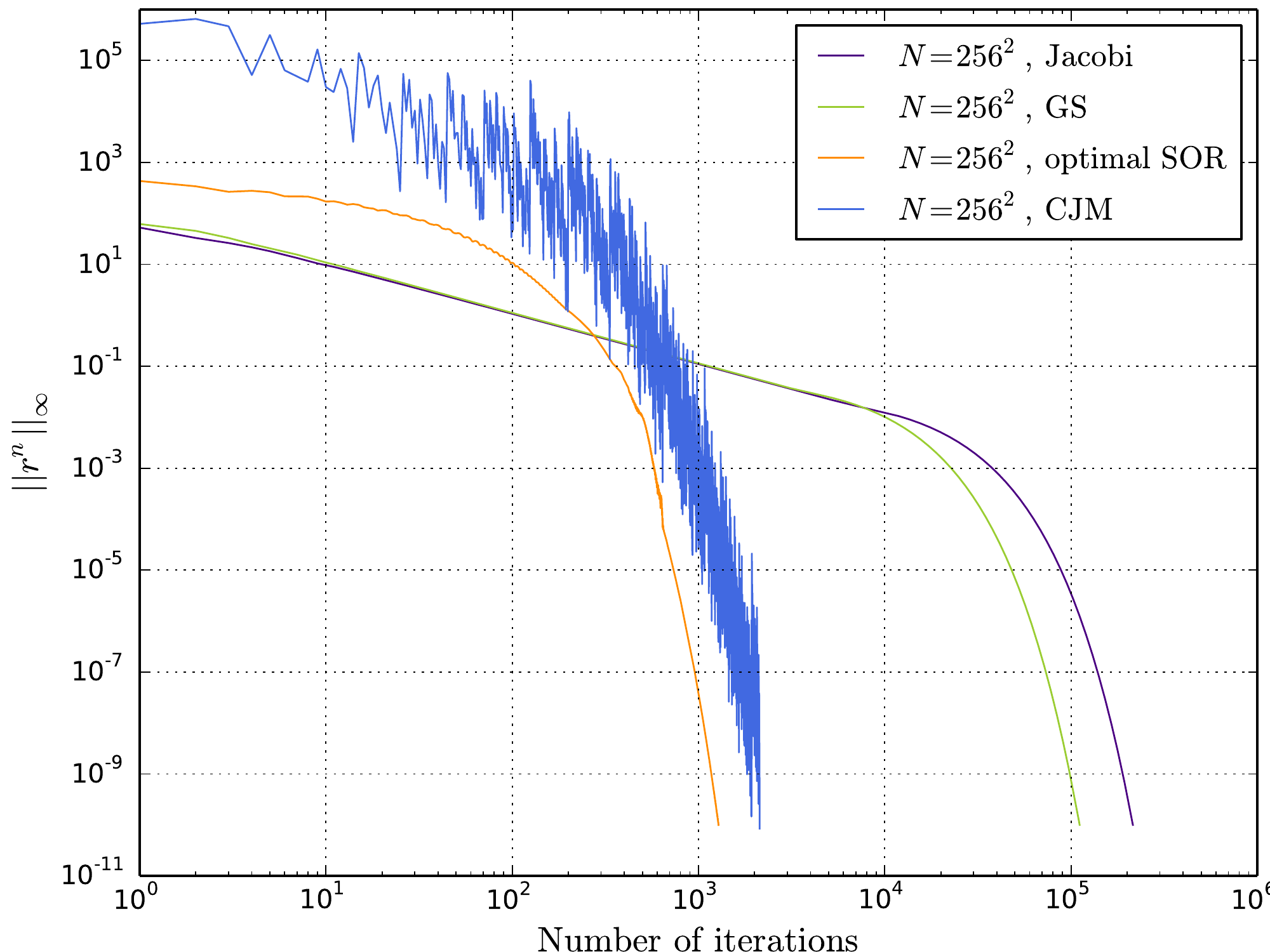} 
\caption{Evolution of the residual for the solution of the Poisson
  equation~\eqref{eq:Poisson2D} in 2D, with 5-points discrete
  Laplacian (Eq.\,\eqref{eq:5-points}; top left panel), 9-points
  (Eq.\,\eqref{eq:9-points}; right panels) and 17-points Laplacian
  (Eq.\,\eqref{eq:17-points}; bottom left panel) for different
  iterative methods and for the resolutions indicated in the
  legens. Note that the top and bottom right panels correspond to a
  problem set up with $N_x=N_y=128$ and $N_x=N_y=256$ points,
  respectively.}
\label{fig:poisson9and17}
\end{figure}

Next, we numerically test the performance of the CJM for the two
high-order discretizations of the Laplacian operator we have
discussed above. To do so, we numerically the following problem:
\begin{equation}
\Delta u = -(x^2 + y^2) e^{x y},
\label{eq:Poisson2D}
\end{equation}
in the unit square with appropriate Dirichlet boundary conditions. The
boundaries are specified easily in this case, since there exists an
analytic solution for the problem at hand that we can compute at the
edges of the computational domain. The analytic solution reads
\begin{equation}
u(x,y)=-e^{x y}.
\end{equation}
In Fig.\,\ref{fig:poisson9and17} we show the residual evolution
obtained when solving problem~\eqref{eq:Poisson2D} with different
high-order discretizations of the Laplacian. In the top left panel we
use the classical 5-points discrete approximation for the Laplacian
(Eq.\,\eqref{eq:5-points}). It is evident that our method almost
reaches the performance of the optimal SOR \cite{Leveque:1988}. In
fact, as we prove in \ref{sec:properties}, this optimal weight for the
SOR method coincides with the geometrical mean of the weights obtained
with our optimal scheme. In the right panels we display the evolution
of the residual when solving the same problem but using the 9-point
discretization of the Laplacian proposed by~\cite{Adams:1988}
(Eq.\,\eqref{eq:9-points}). In the top right panel of
Fig.\,\ref{fig:poisson9and17}, we use a mesh with 128 points in each
dimension, while in the bottom right panel we use 256 points per
dimension. In both cases, the performance is comparable with the
optimal SOR whose weight is calculated in~\cite{Adams:1988}. Finally,
the left-bottom panel of Fig.\,\ref{fig:poisson9and17} shows the
number of iterations when solving the same problem, but using a
$64^2$ grid and our 17-points Laplacian (Eq.\,\eqref{eq:17-points}),
with the optimal CJM obtained with the $\kmin$ and $\kmax$ of
Eqs.\,\eqref{eq:kmin17p} and \eqref{eq:kmax17p} (i.e., in the case $a
= 1,\; b=2$, which gives equal weight to all points in the
neighborhood). In this case, the optimal weight of the SOR is unknown,
so we compute the numerical solution for several values of the SOR
weight. Remarkably, the CJM scheme compares fairly well with SOR.

Last but not least, we are interested in the parallel implementation
of these schemes. It is known that in the case of the standard 5-points
discretization of the Laplacian, one needs to implement a red-black
coloring strategy for the efficient parallel implementation of SOR. In
the case of the 9-points discretization of the Laplacian,
\cite{Adams:1988} points out that one needs four colors for a parallel
implementation. Furthermore, the ordering strategy with more than two
colors is not unique. Adams~\cite{Adams:1988} find 72 different
four-color orderings, which lead to different convergence
rates. In contrast, our CJM scheme (as any SRJ scheme) is trivially
parallelizable since there is no need for a coloring strategy and,
consequently, it possesses a unique convergence rate. We find that the
tiny performance difference between the SOR method, applied to problems
where the optimal weight is unknown, and the CJM is outbalanced by the
simplicity in the parallelization of the latter.

\section{Conclusions}

In this work we have obtained the optimal coefficients for the SRJ method to solve linear systems arising in the finite
difference discretization of elliptic problems in the case $P=M$, i.e., using each weight only once per cycle. We have
proven that these are the optimal coefficients for the general case, where we fix $P$ but allow for repetitions of the
coefficients ($P\le M$). Furthermore, we have provided a simple estimate to compute the optimal value of $M$ to reduce
the initial residual by a prescribed factor.

We have tested the performance of the method with two simple examples (in 2 and 3 dimensions), showing that the
analytically derived amplification factors can be obtained in practice. When comparing the optimal $P=M$ set of
coefficients with those in the literature~\cite{Yang2014,Adsuetal15}, our method always gives better results,
i.e., it achieves a larger reduction of the residual for the same number of iterations $M$. Additionally, the new
coefficients can be computed analytically, as a function of $M$, $\kmax$, and $\kmin$, which avoids the numerical
resolution of the minimization problem involved in previous works on the SRJ. The result is a numerical method
that is easy to implement, and where all necessary coefficients can
easily be calculated given the grid size,
boundary conditions and tolerance of the elliptic problem at hand {\em before} the actual iteration procedure is even started.

We have found that following the same philosophy that inspired the development of SRJ methods, the case $P=M$ results
in an iterative method nearly equivalent to the non-stationary Richardson method as implemented by Young
\cite{Young:1953}; namely, where the coefficients $\omega_{n}$ are taken to be the reciprocals of the roots of the
corresponding Chebyshev polynomials in the interval bounding the spectrum of eigenvalues of the matrix ($A$) of the
linear system. Furthermore, inspired by the same ideas as in the original SRJ methods, the actual minimum and maximum
eigenvalues of $A$ do not need to be explicitly computed. Instead, we resort to a (much simpler) von Neumann analysis
of the linear system which yields the values of the $\kmin$ and $\kmax$ that {\em replace} the (larger) values of the
minimum and maximum eigenvalues of $A$. The key to our success in the practical implementation of the Chebyshev-Jacobi
methods stems from a suitable ordering (or scheduling) of the weights $\omega_{n}$ in the algorithm. Though other
orderings have also been shown to work, our choice clearly limits the growth of round-off errors when the number of iterations is
large. This ordering is inherited from the SRJ schemes.

We have also tested the performance of the CJM for more than second
order discretizations of the elliptic Laplacian operator. These cases
are especially involved since the matrix of iteration cannot be
consistently ordered. Thus, Young's theory cannot be employed to find
the value of the optimal weight of a SOR scheme applied to the
resulting problems. For the particular case of the 9-points
discretization of the Laplacian, even though the iteration matrix
cannot be consistently ordered, Adams~\cite{Adams:1988} found the
optimal weight for the corresponding SOR scheme in a rather involved
derivation. Comparing the results for the numerical solution of a
simple Poisson-like problem of the SOR method derived by Adams and the
CJM we obtain here for the same 9-points discretization of the
Laplacian, it is evident that both methods perform quite similarly
(though the optimal SOR scheme is still slightly better). However, the
SOR method requires a multi-coloring parallelization strategy with up
to 72 four-color orderings (each with different performance), when
applied to the 9-points discretizations of the Laplacian operator. The
parallelization strategy is even more intricate when a 17-points
discretization of the Laplacian is used. In contrast, CJM methods are
trivially parallelizable and do not require any multi-coloring
strategy. Thus, we conclude that the slightly smaller performance
difference between the CJM and the SOR method in sequential
applications is easily outbalanced in parallel implementations of the
former method.

\section*{Acknowledgments}
We acknowledge the support from the European Research Council (Starting Independent Researcher Grant CAMAP-259276) and
SAF2013-49284-EXP, as well as the partial support of grants AYA2015-66899-C2-1-P, AYA2013-40979-P and PROMETEO-II-2014-069.

\appendix

\section{Ordering of the weights}
\label{sec:ordering}

As we point out in Sect.~\ref{sec:OptimalCheb}, the ordering of the weights $\omega_{n}$ is the key to avoid the pile up
of roundoff errors. In this appendix, we show that the ordering provided by \cite{Yang2014} for SRJ schemes, and that we
also use for the optimal $P=M$ schemes, differs from the one suggested by other authors.

Lebedev \& Finogenov \cite{Lebedev:1971} provided orderings for the cases in which the number of weights is a power of
2. Translated to our notation, we shall have $M=2^r$, $r=0,1,\ldots$. In such a case, let the ordering of the set
$(\omega_1, \omega_2, \ldots, \omega_M)$ as obtained from Eq.\,(\ref{eq:omegan}), be mapped with the vector of indices
$(1,2,\ldots,M)$. Let us consider an integer permutation of the vector of indices of order $M$, $\Xi_{M} :=(j_1, j_2, \ldots, j_M)$, where $(1\le
j_k \le M, j_i\ne j_k)$, which are constructed according to the following recurrence relation:
\begin{gather}
\Xi_{2^0}=\Xi_1 := (1) \quad \text{and}\quad  \Xi_{2^{r-1}} :=(j_1, j_2, \ldots, j_{2^{r-1}})\\
 \Xi_{2^{r}} = \Xi_M :=(j_1, 2^r+1-j_1, j_2, 2^r+1-j_2,\ldots, j_{2^{r-1}}, 2^r+1-j_{2^{r-1}})
\end{gather}
In particular, we have,
\begin{eqnarray*}
\Xi_2 &=& (1,2),\nonumber\\
\Xi_4 &=& (1,4,2,3),\nonumber\\
\Xi_8 &=& (1,8,4,5,2,7,3,6),\nonumber\\
\Xi_{16} &=& (1,16,8,9,4,13,5,12,2,15,7,10,3,14,6,11).\nonumber\\
\end{eqnarray*}
In contrast, we can obtain different SRJ schemes, and correspondingly, different orderings, for the same number of
weights, because of the later depend on the number of points employed in the discretization (see
Eq.\,\ref{eq:omegan}). Furthermore, the ordering also depends on the tolerance goal, $\sigma$ (which sets the value of $M$;
Eq.\,\ref{eq:Pmax}). Next we list some of the orderings we can obtain for different discretizations (annotated in
parenthesis in the form $N_x\times N_y$) and values of $\sigma$:
\begin{eqnarray*}
\Xi_2^{\rm SRJ} &=& (1,2),\nonumber\\
\Xi_4^{\rm SRJ} &=& (1,4,3,2), \nonumber\\
\Xi_8^{\rm SRJ} &=& (1,8,5,2,3,7,4,6) \text{ for } (4\times4, \sigma=0.01), \nonumber\\ 
                          &  & (1,8,5,3,6,2,7,4) \text{ for } (8\times8, \sigma=0.15), \nonumber\\
\Xi_{16}^{\rm SRJ} & =& (1,15,9,2,12,3,4,13,5,6,7,8,10,11,14,16) \text{ for } (4\times4, \sigma=2\times 10^{-5}),\nonumber\\
                          &  & (1,16,9,6,12,3,14,7,10,4,13,5,15,8,2,11) \text{ for } (8\times8, \sigma =6\times 10^{-3}). 
\end{eqnarray*} 
which obviously differ from the orderings $\Xi_j$ for $j\ge 4$. 

We note that \cite{Nikolaev:1972} provided also orderings for arbitrary values of $M$, which coincide with those of
Lebedev \& Finogenov \cite{Lebedev:1971} when $M$ is a power of 2 (i.e., $M=2^r$). Finally, more recently, Lebedev \&
Finogenov \cite{Lebedev:2002} have extended their previous work to a larger number of cases (e.g., $M=2^r 3^s$) and
applied also to Chebyshev iterative methods. We remark that the SRJ ordering of the weights can be applied to arbitrary
values of $M$.

\section{Properties of the weights}
\label{sec:properties}

In this appendix we show some algebraic properties of the weights of
the CJM. The first one is that the harmonic mean of the weights equals
the average of the maximum and minimum weight numbers:
\begin{theorem}
Let $\omega_i$ be the weights given by Eq.\,(\ref{eq:omegan}). Then it holds that
\begin{equation}
\frac{1}{M}\sum_{i=1}^M \omega_i^{-1} = \frac{\kmax + \kmin}{2}.
\end{equation}
\end{theorem}
Proof: 
\begin{equation}
	\frac{1}{M}\sum_{i=1}^M \omega_i^{-1} = \frac{(\kmax+\kmin)}{2} 
- \frac{(\kmax-\kmin)}{2 M} \sum_{i=1}^M \cos\left(\frac{\pi(i-1/2)}{M}\right). 
\end{equation}
Let $j\in[1,M/2]$. Since
\begin{equation}
\cos\left(\frac{\pi(j-1/2)}{M}\right) = -\cos\left(\frac{\pi((M-j+1)-1/2)}{M}\right)  ,
\end{equation}
all the terms in the summation cancel out, except the central one in case $M$ is odd. In this last case, 
$M=2n+1$, and the only remaining term is $\cos\left(\frac{\pi(n+1/2)}{2n+1}\right) = \cos\left(\frac{\pi}{2}\right) = 0$. In general, the summation reads
\begin{equation}
	\frac{1}{M}\sum_{i=1}^M \omega_i^{-1} = \frac{\kmax+\kmin}{2}.
\end{equation}

Corollari: Since the relation between the weights of the stationary RM and the CJM is $\hat\omega  = \omega d^{-1}$,
where $D={\rm diag}(A)$, having all its elements equal to $d$, and since $\hat\omega=2/(a+b)$, where $a=\min{(\lambda_i)}$ and $b=\max{(\lambda_i)}$, being $\lambda_i$
  the eigenvalues of matrix $A$, it turns out that
\begin{equation}
\frac{2d^{-1}}{\kmax +\kmin} = \frac{2}{a+b} = \hat\omega.
\end{equation}

\begin{theorem}
Let $\omega_i$ be the weights given by Eq.\,(\ref{eq:omegan}). Then it holds that
\begin{equation}
	\lim_{n \to +\infty} \left[\prod_{i=1}^n \omega_i^{-1}\right]^{1/n} 
= \left(\frac{\sqrt{\kmax}+\sqrt{\kmin}}{2}\right)^2.
\end{equation}
\end{theorem}
Proof: 
We have empirically checked that the sequence $\left[\prod_{i=1}^n \omega_i^{-1}\right]^{1/n}$ is decreasing. Then, the sequence converges to a finite limit if and only if any subsequence converges to the same limit. Let us consider the subsequence $n=2^{(p-1)}, p\in\mathbb{N}$.

Let us define $A=(\kmax+\kmin)/2$, $B=(\kmax-\kmin)/2$, $C_1 = (A^2-B^2/2)$ and $D_1 = B^2/2$. Let us define by recurrence:
\begin{equation}
	C_i = C_{i-1}^2 - D_{i-1}^2/2, \;\; D_i = D_{i-1}^2/2, \;\; i \geq 2.
\end{equation}
\begin{eqnarray}
%
\prod_{i=1}^n \omega_i^{-1} &=& \bigg[A-B\cos\frac{\pi}{2^p}\bigg]
 \bigg[A-B\cos\frac{3\pi}{2^p}\bigg] \ldots \nonumber \\
&&\bigg[A-B\cos\frac{(2^p-3)\pi}{2^p}\bigg] \bigg[A-B\cos\frac{(2^p-1)\pi}{2^p}\bigg] \nonumber \\
	&=& \bigg[(A^2-\frac{B^2}{2}) - \frac{B^2}{2} \cos\frac{\pi}{2^{(p-1)}}\bigg]
\bigg[(A^2-\frac{B^2}{2}) - \frac{B^2}{2} \cos\frac{3\pi}{2^{(p-1)}}\bigg]
\ldots \nonumber \\
&& \bigg[(A^2-\frac{B^2}{2}) - \frac{B^2}{2} \cos\frac{(2^{(p-1)}-1)\pi}{2^{(p-1)}}\bigg] \nonumber \\
	&=& \bigg[C_1 - D_1 \, \cos\frac{\pi}{2^{(p-1)}}\bigg] \bigg[C_1 - D_1 \,
        \cos\frac{3\pi}{2^{(p-1)}}\bigg] \ldots \nonumber \\
&&\bigg[C_1 - D_1 \, \cos\frac{(2^{(p-1)}-1)\pi}{2^{(p-1)}}\bigg].
\end{eqnarray}
The structure we obtain is analogous to that of the initial product, with the change 
$\{A,B,2^{(p-1)}\} \to \{C_1,D_1,2^{(p-2)}\}$. By induction and the definition by recurrence of 
$C_i, D_i$, the above expression reads:
\begin{equation}
	\prod_{i=1}^n \omega_i^{-1} = C_{p-2}^2 - \frac{D_{p-2}^2}{2} = C_{p-1}.
\end{equation}
The problem is reduced to obtain an expression for $C_{p-1}$. It is trivial to check by induction that 
\begin{equation}
D_i = \frac{B^{2^i}}{2^{(2^i-1)}}.
\end{equation}
 Therefore, 
\begin{equation}
\lim_{p\to\infty} [D_{p-1}]^{1/2^{(p-1)}} = B/2.
\end{equation}

Let us define $X_i = C_i/D_i$. From the definition by recurrence of  $C_i, D_i$, we find the following 
recurrence relation: 
\begin{eqnarray}
X_{i+1} = 2X_i^2-1,\qquad X_1 = 1+\frac{8\kmax \kmin}{(\kmax-\kmin)^2} > 1, 
\end{eqnarray}
and by 
induction it can be checked that $X_i>1, \forall i\in\mathbb{N}$. This implies that the sequence 
$X_i^{1/2^i}$ is increasing: $X_i^2 > 1 \Leftrightarrow X_{i+1} > X_i^2 \Leftrightarrow 
X_{i+1}^{1/2^{(i+1)}} > X_i^{1/2^i}$. Moreover, this sequence is bounded above:
\begin{eqnarray}
	X_{i+1} = 2X_i^2 - 1 < 2 X_i^2 \Rightarrow X_i < X_1^{2^{(i-1)}} 2^{(2^{(i-1)}-1)} \nonumber \\
	\Rightarrow X_i^{1/2^i} < [ X_1^{2^{i-1}} 2^{(2^{(i-1)}-1)} ]^{1/2^i} < \sqrt{2X_1}.
\end{eqnarray}
Therefore, the sequence $X_i^{1/2^i}$ is convergent. 

Let us define 
\begin{equation}
L=\lim_{i\to\infty} X_i^{1/2^i}.
\end{equation}
It can be checked by induction that
\begin{equation}
	\left[\frac{\sqrt{X_1+1}+\sqrt{X_1-1}}{\sqrt{2}}\right]^{2^i} = X_i + \sqrt{X_i^2-1}.
\end{equation}
The case $i=1$ is trivially satisfied. Let us assume that the equality holds for $i$, and check that 
it also holds for $i+1$:
\begin{eqnarray}
	&&\left[\frac{\sqrt{X_1+1}+\sqrt{X_1-1}}{\sqrt{2}}\right]^{2^{i+1}} = \left(X_i + \sqrt{X_i^2-1}\right)^2 
\nonumber \\ 
	&=& 2X_i^2 - 1 + 2X_i\sqrt{X_i^2-1} = X_{i+1} + \sqrt{X_{i-1}^2-1}.
\end{eqnarray}
From the previous equality already proven,
\begin{equation}
	X_i^{1/2^i} \leq \frac{\sqrt{X_1+1}+\sqrt{X_1-1}}{\sqrt{2}} = \left(X_i + \sqrt{X_i^2-1}\right)^{1/2^i} 
\leq (2 X_i)^{1/2^i} = 2^{1/2^i} X_i^{1/2^i},
\end{equation}
and taking limits for $i\to\infty$ we get 
\begin{equation}
L \leq \frac{\sqrt{X_1+1}+\sqrt{X_1-1}}{\sqrt{2}} \leq L,
\end{equation}
thus, we conclude that 
\begin{equation}
L = \frac{\sqrt{X_1+1}+\sqrt{X_1-1}}{\sqrt{2}}.
\end{equation}
We finally compute the limit of the geometric mean:
\begin{eqnarray}
	&&\lim_{p\to\infty} \left( \prod_{i=1}^n \omega_i^{-1} \right)^{1/n} = \lim_{p\to\infty} (C_{p-1})^{1/2^{(p-1)}} 
= \lim_{p\to\infty} (D_{p-1})^{1/2^{(p-1)}} (X_{p-1})^{1/2^{(p-1)}} \nonumber \\
	&=& \frac{B \, L}{2} = \left(\frac{\sqrt{\kmax}+\sqrt{\kmin}}{2}\right)^2.
\end{eqnarray}


\bibliographystyle{elsarticle-num} 
\bibliography{cheb}

\end{document}